\documentclass{amsart}
\usepackage[latin1]{inputenc}
\usepackage{amssymb}
\usepackage{latexsym}
\usepackage{longtable}

\newtheorem{teor}{Theorem}[section]
\newtheorem{cor}[teor]{Corollary}
\newtheorem{lema}[teor]{Lemma}
\newtheorem{prop}[teor]{Proposition}

\newtheorem{defn}[teor]{Definition}

\newtheorem{remark}[teor]{Remark}

\newcommand{\Q}{\mathbb{Q}}
\newcommand{\Z}{\mathbb{Z}}

\newcommand{\C}{\mathbb{C}}
\newcommand{\I}{\operatorname{I}}
\newcommand{\PP}{\mathbb{P}}
\newcommand{\Gal}{\operatorname{Gal\,}}
\newcommand{\cond}{\operatorname{cond\,}}
\newcommand{\GL}{\operatorname{GL}}

\newcommand{\Jac}{\operatorname{Jac}}

\newcommand{\id}{\operatorname{id}}
\newcommand{\Aut}{\operatorname{Aut}}
\newfont{\gotip}{eufb10 at 12pt}

\newcommand{\ga}{\mbox{\gotip a}}

\newcommand{\cO}{{\mathcal O}}

\newcommand{\cH}{{\mathcal H}}

\newcommand{\cI}{{\mathcal I}}
\newcommand{\cF}{{\mathcal F}}
\newcommand{\cS}{{\mathcal S}}

\newcommand{\om}{{\omega }}
\newcommand{\ra}{{\rightarrow }}
\newcommand{\hra}{{\hookrightarrow }}
\newcommand{\qbar }{{\bar \Q }}
\newcommand{\Norm}{\mathrm{N}}
\newcommand{\SL}{{\mathrm {SL}}}
\newcommand{\CM}{{\mathrm {CM}}}
\newcommand{\n}{{\mathrm{n}}}
\newcommand{\M}{{\mathrm{M}}}

\newcommand{\R}{{\mathbb R}}

\newcommand{\Pic}{\mathrm{Pic}}

\include{thebibliography}

\begin{document}

\title[Non-elliptic Shimura curves of genus one]
{Non-elliptic Shimura curves of genus one}

\author{Josep Gonz\'{a}lez, Victor Rotger}
\footnote{Both authors are supported in part by by DGICYT Grant
BFM2003-06768-C02-02}

\address{Universitat Polit\`{e}cnica de Catalunya,
Departament de Matem\`{a}tica Aplicada IV (EPSEVG), Av.\ Victor Balaguer s/n,
08800 Vilanova i la Geltr\'{u},
Spain.}

\email{josepg@ma4.upc.edu,vrotger@ma4.upc.edu}



\begin{abstract}
We present explicit models for non-elliptic genus one Shimura curves
$X_0(D, N)$ with $\Gamma_0(N)$-level structure arising from an
indefinite quaternion algebra of reduced discriminant $D$, and
Atkin-Lehner quotients of them. In addition, we discuss and extend
Jordan's work \cite[Ch.\,III]{JoPh} on points with complex
multiplication on Shimura curves.
\end{abstract}

\maketitle

\section{Introduction}

Let $D$ be the reduced discriminant of an indefinite quaternion
algebra $B$ over $\Q $ and and let $N\ge 1$ be a positive integer
coprime to $D$. Let $X_0(D, N)/\Q $ be the Shimura curve over $\Q $
of discriminant $D$ and level $N$ attached to an Eichler order of
level $N$ in $B$. When $D=1$ these are the classical modular curves
$X_0(N)$ which have been extensively studied. Throughout this
article, let us assume $D\ne 1$.

It follows from the genus formula for $X_0(D, N)$ that only for $(D,
N)$ with values in the table below, the genus of $X_0(D, N)$ is $1$.

$$
\begin{array}{|c|c|c|c|c|c|}
\hline (D, N) & (14, 1) & (15, 1) & (21, 1) & (33, 1) & (34, 1) \\
\hline (46,1) & (6, 5) & (6, 7)  & (6,13) & (10, 3) & (10, 7)\\
\hline
\end{array}
$$

\centerline{\ \ \ \ \ {\em Shimura curves $X_0(D,N)$ of genus
$1$.}}
\medskip

By a result of Shimura, $X_0(D, N)(\R )=\emptyset $ and hence these
curves are not elliptic curves over $\Q $.

For trivial level structure $N=1$, equations for all these curves
have already been given (cf.\,\cite{El}, \cite{JoPh}, \cite{JoLi},
\cite{Kur}) except for $D=34$. In some of these cases (especially
for $D>10$), the method employed to construct an equation for
these curves led to large ad hoc computations and messy
diophantine equations (cf. for instance \cite{JoPh}, pp.\,57-68
for the curve $X_0(33, 1)$). These computations turn out to be
even less feasible to handle when one attempts to apply the same
ideas to the discriminant $D=34$.

In Section $2$ we present a simple procedure to provide equations
for curves of genus one provided certain initial data is at our
disposal. In Section $3$ we apply these methods to write down
explicit equations for all the above mentioned curves $X_0(D, N)$
of genus $1$. Since the genus of $X_0(D, N)$ is never $0$ nor $2$
when $N>1$, the present work together with \cite{GoRo},
\cite{JoPh} and \cite{Kur} completes the full list of curves
$X_0(D, N)$ of genus $g\le 2$. In particular, we prove that
Kurihara's conjectured equation \cite{Kur2} for $X_0(34,1)$ is
correct.

Moreover, in Section $4$ we show how our procedure allows us also to
compute equations for the seventeen Atkin-Lehner quotients of
$X_0(D,1)$ of genus $1$ which are non-elliptic over $\Q$. Our
methods also  apply to Atkin-Lehner quotients of Shimura curves
$X_0(D, N)$ with nontrivial level structure $N$, but we do not
include these computations here for the sake of brevity.

As in \cite{GoRo}, \cite{JoPh}, \cite{Kur} we make a crucial use of
the diophantine properties of Shimura curves: their points of
complex multiplication and the class fields generated by them, the
group of Atkin-Lehner involutions acting on $X_0(D, N)$ and their
fixed points, and Cerednik-Drinfeld's description of the special
fibres of $X_0(D, N)$ at primes $p\mid D$ of bad reduction.

However, our approach differs from the previous works in that we
take advantage of an explicit description which goes back to Cassels
of the $\Q$-soluble $\Q$-equivalence classes of an elliptic curve
over $\Q$, and has been made explicit by Cremona and Stoll
\cite{Cremona}, \cite{Stoll}.

In \cite[Ch.\,III]{JoPh} Jordan proves fundamental statements on
complex multiplication points on Shimura curves with trivial level
structure attached to maximal orders of imaginary quadratic
fields. Since in this note we work in the more general setting of
$\Gamma_0(N)$-{\em level structure} and points with complex
multiplication by {\em non maximal} imaginary quadratic orders
(which arise in a natural way for instance as fixed points of some
Atkin-Lehner involutions), we extend these statements to this more
general context in the {\em appendix} to this note.

Most of Jordan's arguments in \cite[Ch.\,III]{JoPh} extend in a
straightforward way, except for the local behavior at primes
dividing both the level $N$ and the conductor $f$ of the quadratic
order. Indeed, primes $p\mid (N, f)$ deserve a closer analysis, as
the approach given in \cite[Ch.\,III]{JoPh} does not apply
immediately to these (cf.\,specially Lemma \ref{CC}). Due to this
and the fact that Jordan's Ph.\,D.\,Thesis \cite{JoPh} is
unpublished and not easily available, we present these results in
the appendix to this article, including proofs of all them in full
generality.

In the last years there have been interesting explicit and
computational approaches to Shimura curves. As recent
contributions let us mention the works of Baba-Granath
\cite{BaGr}, Bayer \cite{Ba} and Elkies \cite{El}, \cite{El2}.
Some of our results may be regarded as progress towards the open
problems posed in \cite{El}.

It is a pleasure to thank Anatoli Segura for a careful reading of
previous drafts of this article and valuable comments on it.

\section{Explicit models for genus one double coverings of
$\PP^1$}\label{models}

Let $C$ be a (projective, nonsingular) curve of genus one over a
field $K$ of characteristic different from $2$. Let us denote by
$\cI(C)$ the set of involutions acting on $C$ over a separable
closure $\overline K$ of $K$, i.e.,
$$\cI(C):=\{v\in \Aut_{\overline K} (C):
v^2=\id\}\,.$$ For $i=0,1$, set
$$
\cI_i(C):=\{v\in \cI(C): C/\langle v\rangle \text{ has genus
$i$}\}\,.
$$

We recall that $\cI_1(C)$ is a group isomorphic to $(\Z/2\,\Z)^2$,
whose elements commute with all involutions of $C$. For a given
$v\in\cI(C)$, set $\cF_v=\{P\in C(\overline{K}):v(P)=P\}$ and let
$K_v$ denote the field extension of $K$ obtained by adjoining the
coordinates of all $P\in \cF_v$. For $v\in\cI(C)$, $v\ne
\mathrm{id}$, it is well known that $\cF_v\neq \emptyset$ if and
only if $v\in\cI_0(C)$ and in this case $\mid \cF_v\mid=4$.
Moreover, for any two different involutions $u, v\in \cI_0$ we have
$\cF_{v}\cap\cF_{u}=\emptyset$, and $u, v$ commute if and only if
$u\cdot v\in\cI_1(C)$.

The following result is well known.

\begin{lema}\label{basic}
The following conditions are equivalent:
\begin{itemize}
\item[(i)] There exists an involution $w\in \Aut_{K} (C)$ such that
$C/\langle w\rangle\stackrel{K}{\simeq }\PP^1_{K}$.

\item[(ii)] There exists $P\in C(\overline K)$ such that $[K(P):K]\leq
2$.

\item[(iii)] There exist $x, y\in K(C)$ and a polynomial $f[X]\in
K[X]$ of degree $3$ or $4$ such that $y^2=f(x)$ and $K(C)=K(x,y)$.
\end{itemize}

\end{lema}

Assuming (ii), we quickly describe these equivalences. For a point
$P$ as in (ii) take $\sigma\in\Gal(\overline{K}/K)$ such that the
divisor $D=(P)+({}^{\sigma}P)$ is defined over $K$. By
Riemann-Roch's Theorem, there exists a nonconstant function $x\in
K(C)$ such that $\operatorname{div}\, x\geq -D$. Since the field
extension $K(C)/K(x)$ has degree $2$, the nontrivial involution of
$K(C)$ over $K(x)$ acts on $H^0(C,\Omega^1)$ as multiplication by
$-1$. The functions $x$ and $y=dx/\omega$, where $\omega$ is a
nonzero regular differential of $C$ defined over $K$, satisfy the
conditions stated in (iii) and $w$ acts on $K(C)$ by sending $(x,y)$
to $(x,-y)$ and $P$ to ${}^{\sigma}P$. The polynomial $f(X)$ has
degree $3$ or $4$ depending on whether $P\in C(K)$ or not. Finally,
for an involution $w$ as in (i),  any  $P \in C(\overline K)$ which
projects to a point in $C/\langle w\rangle(\Q)$ satisfies (ii).  \vskip 0.2
cm

Attached to an equation as in (iii) of the above lemma, there are
two invariants $I$ and $J$ defined by:
$$
I=12\,a_4\,a_0-3\,a_3\,a_1+a_2^2\quad\text{and}\quad
J=72\,a_4\,a_2\,a_0+9\,a_3\,a_2\,a_1-27\,a_4\,a_1^2-27\,a_3^2\,a_0-2\,a_2^3,
$$
where $f(x)=\sum_{i=0}^4 a_i\,x^i$ (cf. \cite{Cremona}). With this
notation, the elliptic curve $E/K$ given by the equation
$$v^2=u^3-27\,I\,u-27\, J$$ is isomorphic over $K$ to the Jacobian
$\Jac (C)$ of $C$.

 From now on, we also  assume that $\mathrm{char }(K)\ne 2, 3$.
For a curve $C/K$ of genus one satisfying the conditions of Lemma
\ref{basic}, let $\pi: C\rightarrow C/\langle w\rangle \simeq
\PP^1_K$ denote the natural projection. The aim of this section is
to present two methods in order to find an equation describing
$C/K$.

\vskip 0.2 cm

\subsection{First method}

In this subsection we describe an approach in order to find an
equation for $C/K$, provided one knows (or can compute) the
following {\em initial data}:

\begin{enumerate}
\item[$(i_1)$] An element $d\in K^*\setminus K^{* 2}$ such that $\pi (C(K(\sqrt{d})))\cap
C/\langle w\rangle(K)\neq \emptyset$.

\item[$(i_2)$] An equation $y^2=x^3+A\, x+B,\,\,
A, B\in K$ for the elliptic curve $E:=\Jac C/K$,

\item[$(i_3)$] The field $K_w$.

\end{enumerate}

For any two polynomials $f_1(x)$, $f_2(x)\in K[x]$ of degree $3$ or
$4$ without double roots, we say that the equations
$eq_1:y^2=f_1(x)$ and $eq_2: y^2=f_2(x)$ are equivalent over $K$ if
there exist
$$\left(\begin{array}{cc} \alpha& \beta\\ \gamma
&\delta\end{array}\right)\in \GL(2,K)\quad \text{and}\quad \lambda
\in K^{*}$$ such that
$$
f_2(x)=\lambda^2 f_1\left (\frac{\alpha \,x+\beta}{\gamma \,
x+\delta}\right)(\gamma \, x+\delta)^4\,.
$$
Let $C_i$ be the curve given by the equation $eq_i$, $i\leq 2$.  If
$eq_1$ and $eq_2$ are equivalent over $K$, then
$C_1\stackrel{K}{\simeq}C_2$ and the splitting fields of $f_1$ and
$f_2$ are equal. However, the converse is not true.

It is known that for a given elliptic curve $E$ over $K$, there is a
one-to-one correspondence between the set $E(K)/2\, E(K)$ and the
set of $K$-equivalence classes of equations $y^2=f(x)$, $\deg (f)=3$
or $4$, such that the corresponding curves given by these ones are
isomorphic to $E$ over $K$ (cf.\,\cite{Cremona}).

Let $C/K$ be a non-elliptic curve of genus one satisfying the
conditions of Lemma \ref{basic}. Assume further that we know the
initial data $(i_1), (i_2), (i_3)$.

Since $C$ admits an equation of the form $y^2=d\, f(x)$, where
$f\in K[x]$ is monic of degree $4$ and the action of $w$ is given
by $(x,y)\mapsto (x,-y)$, we propose the following strategy in
order to determine $f$.

Consider the twisted elliptic curve $E_d: y^2=x^3+ A\, d^2 \,x+ B\,
d^3$ of $E$ and determine a set $\cS=\{\infty,P_1=(x_1,y_1),\dots
,P_r=(x_r,y_r)\}\subseteq E_d(K)$ of representative elements of
$E_d(K)/2 E_d(K)$. The equations  $y^2=f_i(x)$, $0\leq i \leq r$,
where
$$f_i(x)=
 \left\{
 \begin{array}{ll}
x^3+ A\,d^2\, x+ B\,d^3 & \text{if $i=0$,}\\
 x^4-6 x_i \,x^2+ 8 y_i \,x -3 x_i^2 - 4  A \,d^2 & \text{if $1\leq i\leq r$,}
\end{array}
\right.
$$
exhaust all $K$-equivalence classes of equations  attached  to $E_d$
(cf.\,\cite[Proposition 2.2]{Stoll}). Therefore, $C$ must be
isomorphic over $K$ to a curve given by one of the equations
$y^2=d\,f_i(x)$ for $1\leq i\leq r$, since the equation obtained
from $\infty$, $y^2=d\,f_{0}(x)$, corresponds to an elliptic curve
over $K$.

Any diophantine information about $C$ at our disposal may serve to
pick the correct equation for our curve. This may be for instance
the case if the field $K_w$ agrees with exactly one of the
splitting fields of the polynomials $f_i$. As we show in Sections
3 and 4, this approach always succeeds for all Shimura curves and
Shimura curve quotients of genus one that we deal with. Of course,
in general there is no reason to expect that the initial data
$(i_1), (i_2), (i_3)$ suffices to determine $C$.

\subsection{Second method}

Let $C/K$ be a curve of genus one as in Lemma \ref{basic}. In
order to describe a second method to provide an explicit model for
$C$ under additional assumptions, we need the following result.

\begin{prop}\label{model}
Let $C/K$ be a curve of genus one together with $w,u\in \cI_0(C)$
defined over $K$ such that $C/\langle w\rangle \simeq \PP^1_K$ and
that $u\cdot w\in \cI_1(C)$.

Let $L=K(\sqrt{d})$ be a quadratic extension of $K$ and $\sigma
\in \Gal (L/K)$, $\sigma \ne 1$, such that there exists a point
$P\in C(L)\backslash C(K)$ with $w(P)={}^{\sigma}P$. We have that

\begin{enumerate}
\item If $P\in \cF_u$, there exist $x,y\in K(C)$ such that $K(C)=
K(x, y)$, $w(x,y)=(x,-y)$, $u(x, y)=(-x,y)$ and
$$
y^2= d(x^4+ b\, x^2+ c)\,, \quad b,c\in K \,.
$$
Moreover, $ Y^2 = d(X^3 + b X^2+c\,X)$ is an equation for $C/\langle
u\cdot w\rangle$.

\item  If $P\not \in \cF_u$, there exist $x,y\in K(C)$  such that $K(C)=K(x, y)$, $w(x,y)=(x,-y)$,
$u(x,y)=(\varepsilon/x,\varepsilon y/x^2)$ for some
$\varepsilon\in K^*$ and
$$y^2=
d(x^4+ bx^3+ c x^2+b \,\varepsilon x+\varepsilon^2 )\,,\quad
b,c\in K.
$$
Moreover, $Y^2 =d ( X^2-4\varepsilon )(X^2+b\,X+c-2\varepsilon)$
is an equation for $C/\langle u\cdot w\rangle $.
\end{enumerate}
\end{prop}

\noindent {\bf Proof.} Let us first assume that $P\in\cF_u$. Then $\cF_u=\{P,{}^{\sigma}P, Q, Q'\}$ for some
points $Q$, $Q'\in C(\overline K)$ such that the divisor $D=(Q)+(Q')-(P)-({}^{\sigma}P)$ is a $K$-rational
divisor invariant under $w$. Thus, there exist $x,y\in K(C)$ such that $K(C)=K(x,y)$, $$\operatorname{div}\, x=D
\quad \text{and}\quad y^2=a f(x)\,, $$ where $a\in K^{*}$ and $f(X)\in K[X]$ is a monic polynomial of degree $4$
without double roots. Since the value of the function $(y/x^2)^2$ at $P$ and ${}^{\sigma}P$ is $a$, it follows
that $a=d\, a_0^2$ for a certain $a_0\in K^{*}$. Switching $y$ by $a_0\,y$, we can assume that $y^2=d\,f(x)$.
Since $D$ is also invariant under $u$, it follows that $u$ maps $x$ to either $x$ or $-x$. Since $u\in\cI_0$ and
$u\neq w$, we deduce that $u$ acts on $C$ by mapping $(x,y)$ to $(-x,y)$ and that $f$ is an even polynomial. The
function field of $C/\langle u\cdot w\rangle$ is generated by the functions $X= x^2$ and $Y= x\, y$, which
clearly satisfy the equation claimed in our statement.

Assume now that $P\notin\cF_u$. As before, the divisor
$D=(u(P))+(u({}^{\sigma}P))-(P)-({}^{\sigma}P)$ is rational over
$K$ and invariant under $w$. Since $w(P)={}^{\sigma}P$ and $u\cdot
w$ has not fixed points, it follows that $u(P)\neq {}^{\sigma}P$.
Therefore, there exist $x,y\in K(C)$ such that
$$K(C)=K(x,y)\,,\quad w(x,y)=(x,-y)\,,\quad\operatorname{div}\, x=D
\quad\text{and } \quad y^2= d\,f(x)\,, $$ where $f(X)\in K[X]$ is
a monic polynomial of degree $4$. Since $u(D)=-D$, $u$ maps $x$ to
$\varepsilon/x$ and $y$ to $\varepsilon y/x^2$ for some
$\varepsilon\in K^{*}$. Thus $f(X)=X^4+ b\, X^3+ c\,
X^2+\varepsilon\, b\, X+\varepsilon^2$ for some $b, c\in K$ and
the function field of $C/\langle u\cdot w\rangle$ is generated by
the functions $X= x+\varepsilon/x$ and $Y= y(1-\varepsilon/x^2)$,
which again satisfy the equation claimed in our statement. \hfill
$\Box$

\begin{remark}
Under the assumptions of Proposition \ref{model}, there exists (in
both cases) a point $P\in \Jac ( C) [2](K)$ such that $\Jac
(C)/\langle P\rangle=\Jac (C/\langle u\cdot w\rangle)$. In
particular, it holds that $\mid \Aut_K(C)\cap \cI_1\mid$ divides
$\mid \Jac ( C) [2](K)\mid$. If $K$ is a number field, this
implies that both Jacobians have the same conductor over $K$.
\end{remark}

\begin{remark} Assume that $C(K)=\emptyset $. In case {\rm (1)},
$C/\langle u\cdot w\rangle$ is an elliptic curve over $K$. In case
{\rm(2)}, if in addition there exists a point $P\in\cF_u$ such that
$\pi(P)\in\PP_K^1$, then $C$ admits a model as in {\rm (1)} for a
suitable choice of $d$; otherwise, $\varepsilon\notin K^2 $,
$K(\sqrt{\epsilon})$ is a subfield of  $K(\cF_u)$ and $C/\langle
u\cdot w\rangle$ might not be an elliptic curve over $K$.
\end{remark}

In the particular case that $C$ belongs to case $(1)$ of Proposition
\ref{model}, the following result describes an easier and better
procedure to find an equation for $C$, provided one knows the
initial data $(i_1)$, $(i_2)$ as above and

\begin{enumerate}
\item[$(i_4)$] An  equation $V^2= U^3+A'\, U+B'$, $A',B'\in K$,  for the elliptic curve $C/\langle u\cdot
w\rangle $.
\end{enumerate}

\begin{prop}\label{criterion}
Let $C/K$ be as in {\rm (1)} of Proposition \ref{model}. Assume
that $V^2=U^3+A'\, U+B'$, $A', B'\in K$, is an affine equation
over $K$ for the elliptic curve $C'=C/\langle u\cdot w\rangle$.
Then,
$$y^2=d\, x^4+ 3 u_0\, x^2+ (A'+ 3 u_0^2)/d$$
is an equation for $C$, where $u_0\in K$  is the root of the
polynomial  $U^3+A'\, U+B'$  such that $C'/\langle (u_0,0)\rangle$
is isomorphic to $\Jac C$ over $K$.
\end{prop}

\noindent {\bf Proof.} By Proposition \ref{model}, $ Y^2= d\,(
X^3+b\, X^2+ c\, X)$ is an equation for $C'/K$ and it is
isomorphic over $K$ to the elliptic curve
$$V^2=(U-\frac{b\,d}{3}) ( U^2+\frac{b\, d}{3}\,
U-d^2\frac{2\,b^2-9\,c}{9} )= U^3-\frac{d^2(b^2 -
3\,c)}{3}U+\frac{b\,d^3(2\,b^2 - 9\,c)}{27}\,.$$ Computing the $I$
and $J$ invariants attached to the equation $y^2=d(x^4+ b \,x^2+c)$,
it can be checked that the quotient curve $C'/\langle
(b\,d/3,0)\rangle$ is isomorphic over $K$ to the elliptic curve
$v^2=u^3-27\,I\,u-27 J$, which is isomorphic over $K$ to $\Jac C/K$.

Since the $K$-equivalence class of the equation $y^2=d\,x^4+ 3 u_0\,
x^2+ (A'+ 3 u_0^2)/d$ does not depend on the chosen equation
$$V^2=U^3+ A'\, U+B'=\left(U-u_0\right)\left( U^2+u_0\, U+ (A'+u_0^2)\right)$$
for $C'/K$, we can assume the following equalities $$\frac{d\,
b}{3}=u_0\,,\quad -d^2\frac{2\,b^2-9\,c}{9}=A'+ u_0^2\,.$$ It
follows that  $d\, b=3u_0$ and $d\, c=(A'+ 3\, u_0^2)/d$. \hfill
$\Box$

\section{Genus one Shimura curves}\label{g1}

Let $D=p_1\cdots p_{2r}$, $r>0$, be the product of an even number of
distinct prime numbers and let $N\ge 1$, $(D, N)=1$ be an integer.
Let $X_0(D,N)/\Q $ be the canonical model over $\Q $ of the Shimura
curve of discriminant $D$ and level $N$. Let $W_{D, N} =\{ \om_m :
m\mid D\cdot N, (m, D\cdot N/m)=1\} \simeq (\Z /2\Z )^{\sharp \{
p\mid D\cdot N\}}$ be the group of Atkin-Lehner involutions on
$X_0(D, N)$. All these involutions are defined over $\Q $ and we let
$\pi_m: X_0(D, N) \ra  X_0(D, N)/\langle \om_m\rangle $ the natural
projection (cf.\,the Appendix for more details). We shall also let
$K_m=K_{w_m}$ denote the field extension over $K$ obtained by
adjoining the coordinates of the fixed points of $w_m$.

The aim of this section is to provide equations for these curves
when their genera are $1$. Note that these are never elliptic curves
over $\Q $ because they fail to have real points (cf. \cite{Sh2}).

\begin{lema}\label{lema} The Shimura curve $X_0(D, N)$ has genus one
exactly for the  following va\-lues of $(D,N)$: $(14,1),\,(15,1),
\,(21,1),\, (33,1), (34,1), \,(46,1),\,(6,5),\, (6,7),(6,13),\\
(10,3),(10,7)$.
\end{lema}

\noindent{\bf Proof.}  It readily follows from a close inspection
to the genus formula for $X_0(D, N)$ given in Proposition
\ref{genus}. \hfill $\Box$

In order to apply the methods indicated in the previous section,
let us mention what are the key ingredients we use about (genus
one) Shimura curves:

\begin{itemize}
\item The determination of the isogeny class of the elliptic curve $\Jac (X_0(D, N))$ over $\Q$ can be carried
out from Ribet's isogeny theorem (cf. Theorem \ref{Ribet}). In all cases of Lemma \ref{lema}, it turns out that
$\Jac (X_0(D, N))$ lies in the single isogeny class of conductor $D\cdot N$.

\item When the genera of $X_0(D, N)$ and $X_0(D, N)/\langle \om_m\rangle $
are $1$, their Jacobians are isogenous over $\Q $ and their
$\Q$-isomorphism classes are determined by the computation of the
Kodaira symbols of the reduction of both curves at primes $p\mid D$
by using Cerednik-Drinfeld's Theory (\cite[Section 1.7]{BD},
\cite{JoLi}, \cite{Kur}), combined with Table 1 of \cite{Creftp}.
David Kohel's {\em Brandt modules} package implemented in Magma
\cite{magma} is very practical to determine Cerednik-Drinfeld's dual
graphs of $X_0(D, N)$ at primes $p\mid D$.

\item For all curves in Lemma \ref{lema}, there exists an
imaginary quadratic field $K=\Q (\sqrt{d})$ of class number $1$ and
a point $P \in X_0(D, N)(K)$ such that $\pi_{D\cdot N}(P) \in
X_0(D,N)/\langle \om_{D\cdot N}\rangle (\Q )$. The explicit
computation of $d$ follows from Corollary \ref{AL}. In particular, it turns out that $$X_0(D,N)/\langle
\om_{D\cdot N}\rangle \simeq \PP^1_{\Q }$$ for all these curves,
since their genera are always $0$.

\item For every $m\vert D\cdot N$, the number field  $K_m$  can be determined
by using  Proposition \ref{FP}, Theorem \ref{mainCM} and Remark \ref{uniqueness}. These
numbers fields are displayed in the next Lemma.

\end{itemize}

\begin{lema}\label{Km} The set of Atkin-Lehner involutions $\om_m\in \cI_0(X_0(D, N))$ for
the Shimura curves of Lemma \ref{lema} just as the field of
definition of their fixed points are collected in the following
tables:

$$\begin{array}{lll}
(14, 1), \{\om_{14}, \om_2\} \quad & (15, 1), \{ \om_{15}, \om_3\} \quad &
(21, 1), \{ \om_{21},\om_3\} \\[2 pt]
\hline K_{14}= \Q(\sqrt{-1\pm \sqrt{-7}}^{\phantom{c}}) \quad &
K_{15}=\Q(\sqrt{-3})&K_{21}=\Q(\sqrt{-3},\sqrt{-7})
\\[3 pt] K_2= \Q(\sqrt{-1},\sqrt{-2})   &   K_3= \Q(\sqrt{-3}) \quad &K_7= \Q(\sqrt{-7})     \\
\hline
\end{array}
$$

$$\begin{array}{lll}
(33, 1), \{\om_{33}, \om_3\}& (34, 1), \{\om_{34}, \om_{17}\} & (46, 1), \{\om_{46}, \om_2\}\\
\hline K_{33}=\Q(\sqrt{-3},\sqrt{-11}) &K_{34}=\Q(\sqrt{3 \pm
2\sqrt{-2}}^{\phantom{c}})&
K_{46}=\Q(\sqrt{-3\pm \sqrt{-23}})\\
\hline K_3=\Q(\sqrt{-3}) &
K_{17}=\Q(\sqrt{-1\pm 4 \sqrt{-1}}^{\phantom{c}}) & K_{2}= \Q(\sqrt{-1},\sqrt{-2})\\
\hline
\end{array}
$$

$$\begin{array}{lll}
(6,5), \{\om_{30}, \om_2, \om_6, \om_{10}\} & (6,7), \{ \om_{42}, \om_3, \om_6, \om_{21}\} & (6,13),
\{\om_{78}, \om_2, \om_3, \om_{13}\}\\
\hline K_{30}=\Q(\sqrt{-3}, \sqrt{5}^{\phantom{c}}) &
K_{42}=\Q(\sqrt{-2},\sqrt{-3})
& K_{78}=\Q(\sqrt{-3},\sqrt{13})\\
\hline
K_2=\Q(\sqrt{-1})& K_3=\Q(\sqrt{-3})&K_2=\Q(\sqrt{-1}) \\
\hline
K_{6}= \Q(\sqrt{2},\sqrt{-3}^{\phantom{c}})& K_{6}=\Q(\sqrt{2},\sqrt{-3})  &  K_{3}=\Q(\sqrt{-3})\\
\hline K_{10}=
\Q(\sqrt{-2},\sqrt{5}^{\phantom{c}})&K_{21}=\Q(\sqrt{3},\sqrt{-3}) &
K_{13}=\Q(\sqrt{13},\sqrt{-1}) \\
\hline
\end{array}
$$

$$
\begin{array}{ll}
(10,3), \{\om_{30}, \om_2, \om_3, \om_5\}& (10,7), \{ \om_{70}, \om_{5}, \om_{10}, \om_{35}\}\\
\hline
K_{30}=\Q(\sqrt{-6}, \sqrt{10}) & K_{70}=\Q(\sqrt{-14}, \sqrt{10}) \\
\hline
K_{2}=\Q(\sqrt{-2})& K_{5}=\Q(\sqrt{-1},\sqrt{5}^{\phantom{c}})\\
\hline
K_{3}= \Q(\sqrt{-3})& K_{10}=\Q(\sqrt{-2},\sqrt{5}^{\phantom{c}})\\
\hline
K_{5}=\Q(\sqrt{-1},\sqrt{5}^{\phantom{c}})& K_{35}=\Q(\sqrt{5},\sqrt{-7})\\
\hline
\end{array}
$$
\end{lema}

\begin{remark}
When $K_2=\Q(\sqrt{-1},\sqrt{-2})$, two fixed points of $\om_2$ are rational over $\Q(\sqrt{-1})$ while the
coordinates of the other two lie in $\Q(\sqrt{-2})$. In the remaining cases, the Galois closure of the field of
definition of every $P\in \cF_{\om_m}$ is equal to $K_{m}$. The notation $\Q(\sqrt{a\pm\sqrt b\,}\,)$ means that
the field of definition of each $P\in \cF_{\om_m}$ is either $\Q(\sqrt{a+\sqrt b\,}\,)$ or $\Q(\sqrt{a-\sqrt
b\,}\,)$.
\end{remark}

\begin{teor}\label{equ} Equations for the eleven Shimura curves of genus one and
the action of their Atkin-Lehner involutions are collected in the
following tables:

$$
\begin{array}{|c|ccl|}
\hline
(D,N)& y^2&=&\phantom{cccccccccc} f(x) \\
\hline (14,1)&   y^2 &=& - x^4+13\,x^2 -128  \\
 \hline
(15,1)&  y^2 &=  &-3\,x^4 -82\,x^2 -27 \\
 \hline
(21,1) &   y^2&=&-7\,x^4+94\, x^2-343  \\
\hline
(33,1)& y^2 &= &-3\,x^4-10\, x^2-243  \\
 \hline
(34,1) &  y^2&= &- 3\,x^4 +26\, x^3- 53\, x^2-26\, x-3 \\
\hline
(46,1) &  y^2&=&-x^4+45 \,x^2-512 \\
  \hline
(6,5)&   y^2 &=&  -x^4+61\,x^2 -1024\\
\hline
(6,7)&  y^2 &=  & -3\, x^4 - 34\, x^2 - 2187   \\
\hline (6,13) &   y^2&=&-x^4-115\, x^2- 4096 \\
\hline
(10,3)& y^2 &= &-2\, x^4-11\, x^2-32 \\
\hline (10,7) &  y^2&=&-27\, x^4- 40\, x^3+6\, x^2+40\, x- 27\\
\hline
\end{array}
$$
\centerline{\  {\em Table 1. Equations for Shimura curves of genus
one.}}
\medskip

The involution $\om_{D\cdot N}$ maps $(x,y)$ to $(x,-y)$. The action of the remaining Atkin-Lehner involutions
in $\cI_0(X_0(D, N))$ which together with $\om_{D\cdot N}$ generate $W_{D, N}$ are:
$$\begin{array}{cccccc} 14 &
15 & 21& 33& 34 & 46\\
\hline
\om_2(x,y)& \om_3(x,y) &\om_7(x,y) &\om_3(x,y) &\om_{17}(x,y)& \om_2(x,y) \\
(-x,y)& (-x,y)    & (-x,y)   &  (-x,y)   & (-1/x,-y/x^2)
&  (-x,y) \\
\hline
\end{array}
$$
$$\begin{array}{ccccc}
(6,5) & (6,7)& (6,13)& (10,3) & (10,7)\\
\hline
\om_2(x,y)& \om_3(x,y) &\om_2(x,y) &\om_2(x,y) &\om_{15}(x,y) \\[3pt]
(-x,y)& (-x,y)    & (-x,y)   &  (-x,y)   & (\frac{-1}{x},\frac{-y}{x^2}) \\[3pt]
\hline
 \om_6(x,y)& \om_6(x,y) &\om_3(x,y) &\om_3(x,y) &\om_{10}(x,y) \\[3pt]
(\frac{32}{x}, \frac{32 y}{x^2})& (\frac{-27}{ x},\frac{-27
y}{x^2}) & (\frac{64}{x},\frac{64y}{x^2})   &
(\frac{4}{x},\frac{4 y}{x^2})   & (\frac{2x-1}{x-2},\frac{5y}{(x-2)^2}) \\[3pt]
\hline
\end{array}
$$
\end{teor}

\noindent {\bf Proof.} For $X_0(D, N)$ as in Lemma \ref{lema}, set
$w=\om_{D\cdot N}$ and let $d<0$ be an integer such that there
exists $P\in X_0(D, N)(\Q(\sqrt{d}))$, $\pi_{D N}(P)\in
X_0(D,N)/\langle w\rangle (\Q )$.

Assume first that $N=1$. Then, there exists a single Atkin-Lehner
involution $u\ne w$, $u\in\cI_0(X_0(D, 1))$. Next table collects
Cremona's labels for the elliptic curves $\Jac (X_0(D,1))$ and
$\Jac ( X_0(D,1)/\langle u\cdot w\rangle )$ together with some
possible values for $d$:

$$
\begin{array}{|c|c|c|c|c|c|c|c|}
\hline

D  &  d &\Jac (X_0(D,1))& \Jac (X_0(D,1)/\langle u\cdot w\rangle)
  \\
\hline
14 &  -1, -2&  A2 &  A1 \\
\hline
15 &-3& A1&  A2\\
\hline
21 & -7& A2&   A6 \\
\hline
33  &-3&A1 &   A4 \\
\hline
34 &-3& A3&A4\\
\hline
46 & -1,-2& A2& A1\\
\hline
\end{array}\,.
$$

Comparing the above table with the fields $K_m$ of Lemma \ref{Km},
we obtain that $X_0(D, 1)$ for $D=14$, $15$, $21$, $33$ and $46$
correspond to case (1) of Proposition \ref{model}. Corresponding
equations as collected in Table 1 are automatically obtained by
applying Proposition \ref{criterion}.\footnote{Note that equations
for these five curves were already obtained in \cite{JoPh},
\cite{JoLi} and \cite{Kur}. Some of the models proposed there are
different but correspond to equations which are equivalent over $\Q
$ to ours.}

Let us now consider separately the exceptional case $C=X_0(34,1)$,
to which we apply the first method outlined in Section 2.1. As it
follows from the table above, $C$ can be described by an affine
equation of the form
$$y^2= -3\,f(x)\,,$$
for some monic polynomial $f(x)\in \Q[x]$ of degree $4$. Moreover,
its Jacobian is the elliptic curve $34A_3$ given by the equation
$$E: y^2=x^3-\frac{4945}{3} x -\frac{695374}{27}\,.
$$
Consider the twisted curve of $C$:
$$C_{-3}:y^2=  f(x)\,,$$
which is isomorphic over $\Q$ to the elliptic curve $306B3$:
$$E_{-3}:y^2=x^3-{14835}\, x +{695374}\,.
$$
By means of \cite{magma} we obtain that $E_{-3}(\Q)=\langle P_1,
P_2\rangle \simeq \Z/6\Z\times \Z$, where $P_1=(143,1224)$,
$P_2=(63,104)$ and $\mathrm{ord}(P_1)=6$. Thus, since
$|E_{-3}(\Q)/2E_{-3}(\Q)|=4$, the $\Q$-equivalence class of
$y^2=f(x)$  must agree with one of the following three
$\Q$-equivalence classes attached  to $E_{-3}$:
$$
\begin{array}{|c|c|ccl|}
\hline
\text{case}& P\in  E_{-3}(\Q)/2 E_{-3}(\Q) & y^2&= & \phantom{cccccccc}f(x)\\
 \hline
\text{(i)} & 3 P_1=(71,0) &y^2&=&  x^4- 426\,x^2+44217\\
 \hline
\text{(ii)} &P_2=(63,104)&y^2&=&  x^4-378\,x^2+832\,x +47433\\
 \hline
\text{(iii)}& 3P_1+P_2=(35, -468)&y^2&=&   x^4- 210\,x^2 -3744\,x+5565\\
 \hline
\end{array}
$$

Since (ii) is the single case such that the splitting field of the
polynomial $f(x)$ is $K_w$, we conclude that for $f(x)$ as in (ii),
$y^2=-3\,f(x)$ is an equation for $X_0(34,1)$. The map $(x,y)\mapsto
(6\,x-13,12\,y)$ transforms it into the model proposed in Table $1$.

\vskip0.3 cm

Let us consider now the cases for which $N>1$. In all of them, the
conductor of $\Jac (X_0(D,N))$ is $D\cdot N$ and, moreover, $\Jac
(X_0(D,N))[2](\Q )\simeq (\Z/2\,\Z)^2$ since $(\Z/2\,\Z)^3$ is a
subgroup of $\Aut_{\Q}(X_0(D,N))$.

For $(D, N)= (6,5), (6,7)$, $(6, 13)$ and $(10,3)$ there exists an Atkin-Lehner involution $\om_m$ such that
$K_m=\Q(\sqrt d )$ and $m\neq D\cdot N$. More precisely, these are
$$
\begin{array}{ccccc}
( D ,N) &m &d & \Jac (X_0 (D, N))& \Jac (X_0(D, N)/\langle \om_{D\cdot N/m}\rangle)  \\
\hline
(6,5) & 2 & -1&30A6 &  30A3 \\
\hline
(6,7)& 3 &-3 &   42A3&  42A6\\
\hline
(6,13) &2 &-1 & 78A2&   78A1 \\
\hline (10 ,3)& 2&-2 &30A2 &   30A1 \\
\hline
 \end{array}
$$

Equations as claimed in Table $1$ are immediately obtained by
applying again Proposition \ref{criterion}.

As for $(D, N)=(10,7)$ we proceed similarly as we did for $(34,1)$. We have that $E:=\Jac (X_0(10,7))$ is the
elliptic curve $70A_2$,  because this is the single isomorphism class of conductor $70$ with all its $2$-torsion
points rational over $\Q$ and we take the following equation $ y^2=x^3-283\, x -1482$ for $E$. In this case, we
can take $d=-3$ and we have $E_{-3}:y^2=x^3-2547\,x+ 40014$, which is the elliptic curve $630E2$. In addition,
$E_{-3}(\Q)=\langle P_1, P_2, P_3\rangle \simeq (\Z/2\Z)^2\times \Z$, where $P_1 = {(18, 0})$,
 $P_2 = (39, 0)$ and $P_3 = (-17, 280)$. Hence we must check the
following seven $\Q$-equivalence classes of equations attached  to
$E_{-3}$:
$$
\begin{array}{|c|c|ccl|}
\hline
\text{case}& P\in  E_{-3}(\Q)/2 E_{-3}(\Q) & y^2&= &\phantom{cccccccccc}f(x)\\
 \hline
\text{(i)} &  P_1=(18,0) &y^2&=&  x^4- 108\,x^2+9216\\
 \hline
\text{(ii)} &P_2=(39,0)&y^2&=&  x^4-234\,x^2+5625\\
 \hline
\text{(iii)}& P_1+P_2=(-57, 0)&y^2&=&   x^4+342\,x^2 +441\\
 \hline
 \text{(iv)} & P_3=(-17,280) &y^2&=&  x^4+102\,x^2+2240\,x+9321\\
 \hline
\text{(v)} &P_1+P_3=(63, 360)&y^2&=&  x^4-378\,x^2+2880\,  x -1719\\
 \hline
\text{(vi)}& P_2+P_3=(3, -180)&y^2&=&   x^4- 18\,x^2 -1440\,x+10161\\
 \hline
 \text{(vii)}& P_1+P_2+P_3=(123, -1260)&y^2&=&   x^4- 738\,x^2 -10080\,x-35199\\
 \hline
\end{array}
$$

Only for case (iv) the splitting field of $f(x)$ agrees with
$K_{70}$ given in Lemma \ref{lema}. Thus,
$y^2=-3(x^4+102\,x^2+2240\,x+9321)$ is an equation for
$X_0(10,7)$. The transformation $(x,y)\mapsto ((11 x + 1)/(-x +
1),  36 y/(-x + 1)^2)$ yields the model collected in Table $1$.

Finally, for each of the equations in Table 1, we can compute all
their involutions over $\Q$ which commute with $\om _{D\cdot N}$.
The explicit expression of the Atkin-Lehner involutions acting on
each of these models is determined by comparing these computations
with the fields $K_m$ as in Lemma \ref{Km}. \hfill $\Box$

\begin{remark}
Kurihara conjectured in \cite{Kur2} the following equation for
$X_0(34,1)$:
$$
\left\{\begin{array}{ccc} z^2+ 44\,u^2-68\, u+27&=&0 \,,\\
w^2-(u^2+1)&=&0\,.
\end{array}\right.
$$
It is checked that the map
$$(u,w,z)=\left (-\frac{2 x}{x^2-1},-\frac{1+x^2}{x^2-1}, \frac{9 y}{x^2-1}\right)$$
transforms it into
$$-3 y^2= x^4+ 136/27 x^3+122/27 x^2-136/27 x+1\,.$$
In turn, the transformation $(x,y)\mapsto (\frac{-5\,x+3}{3\, x+5},
\frac{68 y}{9(3\,x+5)^2})$ shows that Kurihara's conjectured curve
coincides with ours collected in Table $1$. See \cite{BFGR} for an
application of the explicit knowledge of an equation for $X_0(34,
1)$.
\end{remark}

\section{Genus one Atkin-Lehner quotients of Shimura curves}

Let $D = p_1\cdots p_{2r}$, $r>0$, be the product of an even number
of distinct prime numbers and let $X_D=X_0(D, 1)$. For a positive
divisor $m\mid D$, $m>1$, let us denote by $X_D^{(m)}$ the
Atkin-Lehner quotient $X_D/\langle \om_m\rangle $. Note that despite
$X_D(\R )=\emptyset $, it might (and does in several cases) happen
that $X_D^{(m)}(\Q )\ne \emptyset $.

\newpage
The complete list of values of $(D, m)$ for which $X_D^{(m)}$ has genus one, together with Weierstrass models
for those $X_D^{(m)}$ which are elliptic curves over $\Q $, can be found in \cite{Ro}\footnote{In Table $2$ of
\cite{Ro}, we claimed that the genus one curves $X_{D}^{(m)}$ for $(D, m)=(35, 7)$, $(51, 3)$ and $(115, 23)$
fail to have rational points over the local fields $\Q_5$, $\Q_{17}$ and $\Q_5$, respectively. This is wrong:
these curves admit rational points everywhere locally. In fact, both three curves are elliptic curves over $\Q
$, since in each case $h(\Q (\sqrt{-D}))=2$ and there exists a point $P\in \CM(\Z [\frac{1+\sqrt{-D}}{2}])$ on
$X_D$ that projects onto a rational point on $X_{D}^{(m)}$. Namely, these three elliptic curves are, in Cremona
notation, $35A1$, $51A2$ and  $115A1$. Finally, let us also note in passing that Table $3$ of \cite{Ro} should
read $(26, \om_{26})$ instead of $(26, \om_{13})$.}. In this section we provide explicit equations for the
remaining genus one curves, that is, those that fail to have rational points over $\Q $.

\begin{lema}\cite{Ro}\label{al1}
The curve $X_D^{(m)}$ is a non-elliptic curve of genus one over $\Q $ exactly for the following values of $(D,
m)$: $\,\,(39, 13),\, (55, 5), (62, 2),\, (69, 3), (77, 11),\\ (85, 17)\,,\, (94, 2)\,, \,\,(178,\, 89)\,,\,\,
(210, \,30)\,,\,\, (210,
\,42)\,, \,\,(210,\, 70)\,,\, (210,\, 105)\,, \,\,(330, \,3), \,\\
(330, \,22), (330, 33), (330, 165), (462, 154)$.
\end{lema}

Let us mention which are the main tools we use in order to apply
the methods exposed in Section \ref{models}:

\begin{itemize}

\item Similarly as in Section \ref{g1}, the isogeny class of the Jacobian of $X_D^{(m)}$ can be
determined by combining Theorem \ref{Ribet} with Table $3$ of
\cite{Creftp}. The isomorphism class is obtained by comparing Table
$1$ of \cite{Creftp} with the Kodaira symbols of the special fibres
at primes $p\mid D$, which can be computed by means of
Cerednik-Drinfeld's Theory.

\item For any $m'\mid D$, the Atkin-Lehner involution $\om_{m'}$ on $X_D$ induces an involution on
$X_D^{(m)}$ that we shall denote $\tilde{\om}_{m'}$. Then
$\cF_{\tilde{\om}_{m'}}=\pi_m(\cF_{\om_{m'}})\cup
\pi_m(\cF_{\om_{m}\cdot \om_{m'}})$. In all cases of Lemma
\ref{al1} it turns out that
$\cF_{\tilde{\om}_{D}}=\pi_m(\cF_{\om_{D}})$. The number field
generated by the coordinates of the fixed points of
$\tilde{\om}_{m'}$ on $X_D^{(m)}$ can be computed by means of
Proposition \ref{FP}, Corollary \ref{AL} and Remark \ref{uniqueness}.

\item For all curves $X_D^{(m)}$ in Lemma \ref{al1}, $X_D^{(m)}/\langle
\tilde{\om}_D\rangle = X_D/\langle \om_m, \om_D\rangle \simeq
\PP^1_{\Q }$. This can be deduced from the following two general
results: Firstly, the set of fixed points of $\om_D$ is nonempty
(cf.\,Proposition \ref{og} and \ref{FP}) and thus $X_D^{(m)}/\langle
\tilde{\om}_D\rangle $ has genus $0$. Secondly, $X_D/\langle
\om_D\rangle $ has points everywhere locally by \cite[Theorem
3.1]{RSY}. By the Hasse principle this implies that
$X_D^{(m)}/\langle \tilde{\om}_D\rangle \simeq \PP^1_{\Q }$.

\item There always exists an imaginary quadratic field $K=\Q (\sqrt{d})$
such that $h(K)=1$ and $\pi (X_D(K))\cap (X_D^{(m)}/\langle
\tilde{\om}_D\rangle (\Q ))\ne \emptyset $, where now we let $\pi:
X_D\rightarrow X_D^{(m)}/\langle \tilde{\om}_D\rangle $ denote the
natural projection of degree $4$.
The computation of $d$ follows from Corollary \ref{AL}.
\end{itemize}

\begin{teor} Equations for the Atkin-Lehner quotients $X_D^{(m)}$ listed in
Lemma \ref{al1} together with Cremona's label for their Jacobians
are collected in the following table
\newpage
$$
\begin{array}{|c|ccl|c|}
\hline (D,m) & y^2 & =& \phantom{ccccccccc}f(x) & \Jac (X_D^{(m)})\\
\hline (39, 13)& y^2 & =& - 7\,x^4 - 24\,x^3- 34\,x^2+ 24\,x -7 & 39A1 \\
\hline
  (55, 5)& y^2 & =&   - 3\,x^4- 2\,x^3- 9\,x^2 + 2\,x -3 & 55A1\\
\hline (62, 2)& y^2 & =&    - x^4- 8\,x^3- 78\,x^2+ 248\,x-961 & 62A3 \\
\hline
 (69, 3)& y^2 & =&  - 3(x^4- 156\, x^2+ 6912) & 69A2\\
\hline (77, 11)& y^2 & =& - 11\,x^4-19\, x^2-16 & 77C2\\
\hline (85, 17)& y^2 & =&  - 3\,x^4+ 10\,x^3- x^2 - 10\,x -3 & 85A1 \\
\hline (94, 2)& y^2 & =& - x^4+9\, x^2-32 & 94A2\\
\hline (178, 89)& y^2 & =&     - 12\,x^4 + 4\,x^3- 19\,x^2- 4\,x -12 & 178B1\\
\hline (210, 30)& y^2 & =&     - 43\,x^4 - 686\,x^3- 2915\,x^2- 1372\,x-172 & 210E5\\
\hline (210,42)& y^2 & =&    - 43\,x^4+ 600\,x^3- 986\,x^2- 600\,x -43 & 210A6\\
\hline(210,70)& y^2 & =&   - 43\,x^4 + 256\,x^3- 130\,x^2- 768\,x -387 & 210C3\\
\hline (210, 105)& y^2 & =&  - 43\,x^4 + 2\,x^3+ 505\,x^2 + 12\,x-1548 & 210B6 \\
\hline(330, 3)& y^2 & =& - 3\,x^4 - 22\,x^3 - 125\,x^2- 66\,x  -27  & 330B3\\
\hline(330, 22)& y^2 & =&  - 3\,x^4+ 1358\,x^2  -177147 & 330C3\\
\hline(330, 33)& y^2 & =& - 3\,x^4+2846 \, x^2-2381643 & 330D2\\
\hline (330, 165)& y^2 & =& - 3\,x^4+614 \, x^2-19683 & 330A2\\
\hline (462, 154)& y^2 & =& - \,x^4-283 \, x^2-16384 & 462B2\\
\hline
\end{array}
$$
\centerline{\ \ \ \ \ {\em Table 2. Equations for non-elliptic
Atkin-Lehner quotients of genus one.}}
\vskip 0.1 cm
where $\tilde{\om}_{D}$ acts as $(x,y)\mapsto (x,-y)$.
\end{teor}
\noindent {\bf Proof.} The equations in Table 2 are obtained by
applying the first procedure described in Section 2. For  every pair
$(D,m)$ we take as $d$ the leading coefficient of the polynomial
$f(x)$ in Table 2.  In all these cases there is a single element in
$E_d(\Q)/2 \, E_d(\Q)$ such that the splitting field of its attached
equivalent class agrees with the field $K_{\tilde{\om}_D}$. We
summarize the computations in the next table:
$$
\begin{array}{|c|c|c|}
\hline (D,m) &(A,B) &  P \\
\hline (39, 13)&
(-217/3, -5510/27) & (-11/3,288)\\
\hline   (55, 5)& (-67, 126)  &(-17,44) \\
 \hline (62, 2)& (-491,-154) &(-9,62),\\
\hline  (69, 3)& (-2235,40534) & (26,0)\\
\hline (77, 11)& (- 2473/3,  227050/27))& (-418/3,0) \\
\hline (85, 17)& (-409/3, -16454/27) &(23,20)\\
\hline (94, 2)& (-155,-714)& (6,0)\\
\hline (178, 89)&(-2137/3, 170170/27)& (-37,-128)\\
\hline (210, 30)&(-5762401/3, -27665272798/27) &(102257/3, 102900) \\
\hline (210,42)& (-129649/3, -90882286/27)  &(46301/3, 1646400)\\
\hline (210,70)& (-50401/3, -22628702/27)  &(9493/3, -4800)\\
\hline (210, 105)& (-1053721/3, -2163135238/27)  &(43433/3, -43904)\\
 \hline(330, 3)& (-12241/3, -2249422/27) &(-129,-2200)\\
 \hline(330, 22)& (-513841/3, -733647278/27)  & (679,0)\\
\hline(330, 33)& (-5864929/3, -22155907934/27)& (1423,0)\\
\hline (330, 165)& (-67849/3, -33554486/27) & (307,0)\\
\hline (462, 154)&  (-276697/3,288510010/27)&(-566/3,0) \\
\hline
\end{array}
$$
\centerline{\ \ \ \ \ {Table 3.}}
\newpage
Here,  $E=\Jac (X_D^{(m)}): y^2=x^3+A\, x+B$  and $P\in E_d(\Q)$
lies in the class of $E_d(\Q)/2 E_d(\Q)$ which provides the single
$\Q$-equivalence class $d\, y^2=f(x)$ isomorphic to $E_d$ such that
the splitting field of $f$ is $K_{\tilde{\om}_D}$. \hfill $\Box$

\begin{remark}
For all these curves, there is at the least one involution $u\in
\cI_0$ defined over $\Q$ commuting with $w$ (though in some cases
$u$ does not arise from any Atkin-Lehner involution on $X_0(D,1)$).
The equations obtained  from the point $P$ in Table 3 have been
replaced in Table 2 by equivalent equations as in (1) or (2) of
Proposition \ref{model} depending on whether there exists an
involution $u$ with a fixed point that projects onto a rational
point on $X_D^{(m)}/\langle \tilde{\om}_D\rangle$ or not.

\end{remark}

\begin{remark}
Kurihara conjecured in \cite{Kur2} equations for the genus three
curves $X_{39}$, $X_{55}$, $X_{62}$, $X_{69}$ and $X_{94}$. For $(D,m)=(39,13),(55,5),(62,2),(69,3), (94,2)$, our theorem above proves that Kurihara's conjectural equations for these genus one quotients are correct.
\end{remark}

\section{Appendix: Shimura curves and their points of complex
multiplication}\label{CM} Let $B$ be an indefinite division
quaternion algebra of discriminant $D=p_1\cdots p_{2 r}$, $D\ge 1$.
For any integer $N\ge 1$ coprime to $D$, Shimura introduced a
projective smooth algebraic curve $X_0(D, N)/\Q $ which can be
described as follows.

Let $\cO_{D, N}$ be an Eichler order of level $N$ in $B$. Let $\n
:B\,\ra \,\Q $ denote the reduced norm on $B$ and let $B^*_+$ be the
subgroup of elements of $B^*$ of positive reduced norm. Let $\cO
_{D, N}^1 = \{ \gamma \in \cO_{D, N}: \n (\gamma ) =1 \}$, which we
regard as a discrete subgroup of $\SL _2(\R )$ through a fixed
isomorphism $\Psi : B\otimes \R \simeq \mbox{M}_2(\R )$. Let $\cH $
denote Poincar\'e's upper half-plane. Then
$$
\cO_{D, N}^1\backslash \cH
$$
is a Riemann surface which is compact unless $D=1$. We let $\Phi :
\cH \rightarrow \cO_{D, N}^1\backslash \cH$ denote the natural
uniformization map.

The following fundamental result is due to Shimura.

\begin{teor}\cite[Main Theorem I]{Sh67}, \cite[Theorem 2.5]{Sh3}
Let $D=p_1\cdots p_{2 r}\ge 1$ be a square-free integer and let $N$
be a positive integer coprime to $D$. There is a projective
algebraic curve $X_0(D,N)/\Q $ such that there exists an open
immersion of Riemann surfaces
$$\cO_{D, N}^1\backslash \cH\hookrightarrow X_0(D,N)(\C ).$$

When $D>1$, this is a biregular isomorphism.

\end{teor}

In the theorem, $X_0(D,N)/\Q $ denotes Shimura's canonical model
over $\Q $ as in \cite[Section 3]{Sh67}. When $D=1$,
$X_0(N):=X_0(1,N)$ stands for the classical elliptic modular curve.

\begin{prop}\label{genus}\cite[p.\,280,
301]{Ogg1}

For $D\ne 1$, the genus of $X_0(D, N)$ is
$$
g=1+\frac{D N}{12}\cdot \prod_{p\mid D}(1-\frac{1}{p})\cdot
\prod_{p\mid N}(1+\frac{1}{p})-\frac{e_3}{3}-\frac{e_4}{4},
$$
where for $k=3, 4$: $$e_k=\prod_{p\mid D}(1-(\frac{-k}{p}))\cdot
\prod_{p\parallel N}(1+(\frac{-k}{p}))\prod_{p^2\mid N}\nu_p(k),
\quad \nu_p(k)=\begin{cases}2 & \text{ if } (\frac{-k}{p})=1 \\
0 & \text{ otherwise.}\end{cases}$$ Here $(\frac{\cdot }{\cdot })$
stands for the Kronecker quadratic symbol.
\end{prop}

\begin{prop}\cite[Proposition 4.4]{Sh2}\label{real}
Let $D>1$, $N\ge 1$, $(D,N)=1$. Then
$$X_0(D,N)(\R )=\emptyset .$$
\end{prop}

A fortiori, curves $X_0(D,N)$ fail to have rational points over $\Q
$ when $D>1$.

Assume for the rest of this appendix that $N$ is square-free. As a
natural subgroup of the group of automorphisms of $X_0(D, N)$ over
$\Q $ there is the Atkin-Lehner group of involutions
$$
W(D,N) = \mathrm{Normalizer}_{B_+^*}(\cO _{D, N}^1)/(\Q ^*\cdot \cO
_{D, N}^1).
$$

Its elements can be labelled as $W(D, N) = \{ \om_m(D,N): m\mid D\cdot N, m>0\}$,
where $\om_m(D,N)\in \cO_{D, N}$ can be taken to be any generator of the only two-sided ideal
of reduced norm $m$ of $\cO_{D, N}$. When no confusion on the choices of $D$
and $N$ can arise, we will simply denote $\om_m=\om_m(D,N)$. The
Atkin-Lehner group is abelian, $W(D,N)\simeq (\Z /2\Z )^{\sharp
{p\mid D N}}$ and $\om_m\cdot \om_n=\om_{m\cdot n/(m,n)^2}$ for any
pair of divisors $n,m\mid D N$.

Let $J_0(D, N)/\Q $ denote the Jacobian variety of $X_0(D,N)$. In
particular, $J_0(N)=J_0(1,N)$ stands for the Jacobian variety of
$X_0(N)$. By the universal property of $J_0(D,N)$, we can regard
$W(D,N)$ as a subgroup of $\Aut_{\Q }(J_0(D,N))$.

For any integer $M\ge 1$ and a positive divisor $d\mid M$, let
$J_0(M)^{d-\mbox{new}}/\Q $ denote the optimal quotient variety of
$J_0(M)$ which is $d$-new with respect to the action of the Hecke
algebra in the sense of \cite[Section 1.7]{BD} and \cite{Rob}. The
action of the group $W(1, M)$ on $J_0(M)$ restricts to a
well-defined action on $J_0(M)^{d-\mbox{new}}$.

\begin{teor}\cite{Ri},\cite[Sections 1.3-1.8]{BD}\label{Ribet} There
exists an isogeny defined over $\Q $
$$
\psi : J_0(D\cdot N)^{D-\mbox{new}} \longrightarrow J_0(D,N)
$$
such that, for each $\om_m(D, N)\in W(D,N)$, we have
$$\psi^*(\om_m(D,N))=(-1)^{\# \{p\mid (D, m)\}}\om_m(1,D\cdot N)\in \Aut_{\Q }(J_0(D\cdot
N)).$$

\end{teor}

Note that when the genus $g(X_0(D,N))=1$ and $D>1$, $X_0(D,N)$ is a
non-elliptic genus one curve, which becomes isomorphic over $\qbar $
to the elliptic curve $J_0(D,N)$. The above result shows that in
this case the conductor of $J_0(D,N)$ is $D\cdot N_0$ for some
positive divisor $N_0\mid N$. It turns out that for all these cases
(cf.\,Lemma \ref{lema}) $\dim J_0(D\cdot N)^{new} = 1$ and thus
$\cond (J_0(D,N)) = D\cdot N$.

\subsection{CM-points and their fields of definition}

Let $K$ be an imaginary quadratic field and let $R\subset K$ be an
order of $K$. Following Eichler, we say that an embedding
$q:R\hookrightarrow \cO_{D,N}$ is {\em optimal} if $q(K)\cap
\cO_{D,N} = q(R)$. Via $\Psi \circ q$, $R\setminus \{ 0\}$ embeds in
$\GL_2^+(\R ) = \{ A\in \GL_2(\R ): \det (A)>0\} $ and there is a
single point $x_q\in \cH $ which is fixed by the action of
$R\setminus \{ 0\} $ on $\cH $. We say that $q$ is {\em normalized}
if for any $a\in K^*$, $\Psi \circ q (a)\cdot
\begin{pmatrix} x_q\\1 \end{pmatrix} = a\cdot\begin{pmatrix} x_q\\1 \end{pmatrix}$.

\begin{defn}
The set $\mathrm{CM}(R)$ of complex multiplication (CM) points by
$R$ on $X_0(D,N)$ is the set $\{ \Phi (x_q)\in X_0(D, N)(\C )\} $,
where $q: R\hookrightarrow \cO_{D, N}$ is any normalized optimal
embedding of $R$ into $\cO_{D, N}$.\end{defn}

As it follows from the definition, a point $x\in \cH $ has complex
multiplication by $R$ if and only if the stabilizer
$$\mathrm{Stab}_{\cO^+_{D, N}}(x) := \{ \alpha \in \cO_{D, N} \cap
B^*_+ : \Psi (\alpha ) \cdot x = x \}$$ is $R\setminus \{ 0\}$.
Indeed, if we let $\Phi (x_q)\in \mathrm{CM}(R)$ for some optimal
embedding $q:R\hra \cO_{D, N}$, we have $\mathrm{Stab}_{\cO^+_{D,
N}}(x_q) = q(R)\setminus \{ 0\}$. Otherwise, if $x\in \cH $ is not a
CM-point, then $\mathrm{Stab}_{\cO^+_{D, N}}(x) = \Z \setminus \{
0\}$.

From now on, we fix the following notation.

\begin{itemize}
\item $p$ denotes an integer prime.

\item $K=\Q(\sqrt{-s})$ is an imaginary quadratic field.

\item $R$ is an order in $K$.

\item $f$ is the conductor of $R$.

\item $(\frac{K}{p})$ is the Kronecker symbol.

\item $(\frac{R}{p})=\begin{cases} (\frac{K}{p}) & \text{ if } p\nmid f \\
1 & \text{ if } p\mid f \end{cases}$ is the Eichler symbol.

\item $\I (R)$ is the group of fractional invertible ideals of $R$.

\item $\mathfrak a,\mathfrak b,\dots $ are elements in $\I (R)$.

\item  $\Norm_{K/\Q}(\ga)=\vert R/\ga\vert$.

\item  $H_R$ is the ring class field of $R$, that is, the
abelian extension of $K$ unramified outside $c$ such that
$\Gal(H_R/K)\simeq \Pic (R)$.

\item $h(R)=[H_R:K]$.

\item  $\sigma _{\mathfrak a}$  is the element in $\Gal(H_R/K)$ attached
to $\mathfrak a$ by the Artin symbol.

\end{itemize}

Attached to the quadratic order $R$, the discriminant $D$ and the
level $N$, let us define
$$
D(R)=\prod _{p\mid D, (\frac{R}{p})=-1} p\,, \quad N(R)=\prod
_{p\mid N, (\frac{R}{p})=1} p\,, \quad N^*(R)=\prod _{p\mid N,
p\nmid f, (\frac{R}{p})=1} p\,,
$$
and $W(R) = \{ \om_m \in W: m\mid D(R) N(R) \}$. Note that
$$\gcd (D(R)\, N^*(R), \text{disc} (R))=1\quad \text{ and}\quad \gcd
(D(R)\, N(R), \text{disc} (R))=\gcd (N, f)\,.$$

\begin{prop}\label{og}\cite[Section 1]{Ogg1}, \cite{JoPh}, \cite[Lemma 2.5]{BD} The set
$\CM(R)$ is nonempty if and only if $\frac{D\, N}{D(R)\, N^*(R)}$
divides $\operatorname{disc} (R)$. Moreover, in this case we have that
$W(R)\times \Gal (H_R/K)$ acts freely and transitively on
$\mathrm{CM}(R)$, and thus

$$\sharp \CM(R)= 2^{\sharp \{ p\mid D(R) N(R)\} }\cdot h(R).$$

\end{prop}

Let $E(R, \cO_{D, N})$ denote the set of normalized optimal
embeddings $q: R\hookrightarrow \cO_{D, N}$. The unit group
$\cO^1_{D, N}$ acts on $E(R, \cO_{D, N})$ by conjugation and there
is a one-to-one correspondence between $\cO^1_{D, N}\backslash E(R,
\cO_{D, N})$ and $\CM(R)$. Proposition \ref{og} above follows from
this correspondence and Eichler's theory on optimal embeddings.

For any prime $p$, let $E_p(R, \cO_{D, N}):=\{ q: R\otimes \Z_p \hra
\cO_{D, N}\otimes \Z_p \}$ be the set of local optimal
embeddings\footnote{There is not a natural notion of {\em
normalized} embeddings into $\cO_{D, N}\otimes \Z_p$. Note also that
the coset $\cO^1_{D, N}\backslash E(R, \cO_{D, N})$ is in one-to-one
correspondence with the set of (non necessarily normalized) optimal
embeddings $q: R\hra \cO_{D, N}$ up to conjugation by $\cO^*_{D,
N}$.} at $p$. The coset $(\cO_{D, N}\otimes \Z_p)^*\backslash \,
E_p(R, \cO_{D, N})$ has cardinality $1$ or $2$, and it is $2$ if and
only if $p\mid D(R)\cdot N(R)$. For any such prime $p$, there is a
natural {\em orientation} map $$o_p: \cO^1_{D, N}\backslash E(R,
\cO_{D, N})\longrightarrow (\cO_{D, N}\otimes \Z_p)^*\backslash \,
E_p(R, \cO_{D, N}) = \{ \pm 1\}.$$ Proposition \ref{og} can be
refined to claim that for any $v\in \{ \pm 1\}^{\sharp \{ p\mid D(R)
N(R)\} }$, the cardinality of the fibre at $v$ of $\prod_{p\mid D(R)
N(R)} o_p$ is $h(R)$. We say that two points $P, P'\in \CM(R)$ lie
in the same {\em branch} if the corresponding normalized optimal
embeddings $q, q'\in E(R, \cO_{D, N})$ are locally equivalent, that
is, have the same orientation at all $p\mid D(R) N(R)$. The set
$\CM(R)$ is the disjoint union of $2^{\sharp \{ p\mid D(R) N(R) \}}$
branches, consisting of $h(R)$ points each (cf.\,\cite[Section
1]{Ogg1}).

More precisely, Galois elements $\sigma \in \Gal(H_R/K)$ preserve all local orientations of points in $\CM(R)$. For $p\mid D(R) N(R)$,
Atkin-Lehner involutions $\om_p$ switch the local orientation at $p$ and preserve the remaining ones (\cite[Lemmas 2.4 and 2.5]{BD}).

Fixed points of Atkin-Lehner involutions acting on Shimura curves
are points of complex multiplication, sometimes by a non-maximal quadratic order:

\begin{prop}\cite[Section 1]{Ogg1}\label{FP} Let $m\mid D\cdot N$,
$m>0$. The set of fixed points of the Atkin-Lehner involution
$\om_m$ acting on $X_0(D, N)$ is
$$
\cF_{\om_m}=\begin{cases} \CM(\Z [\sqrt{-1}])\cup \CM(\Z
[\sqrt{-2}]) & \mbox{ if } m=2 \\ \CM(\Z [\sqrt{-m}])\cup \CM(\Z
[\frac{1+\sqrt{-m}}{2}]) & \mbox{ if } m\equiv 3\mbox{ mod } 4
\\ \CM(\Z [\sqrt{-m}]) & \mbox{ otherwise. }
\end{cases}
$$

\end{prop}

For the rest of this section, let $R\subset K=\Q(\sqrt{-s})$, $s>0$,
be an imaginary quadratic order such that $\CM(R)\ne \emptyset $.

\begin{teor}\label{Shi}\cite[Main Theorem II]{Sh67}, \cite{Sh3}

Let $P\in \operatorname{CM}(R)$ and $\Q (P)$ be the number field
generated by the coordinates of $P$ on $X_0(D, N)$. Then

\begin{enumerate}

\item $H_R=K\cdot \Q (P)$

\item (Shimura's reciprocity law) Let $q: R\hookrightarrow \cO_{D, N}$ be a normalized
optimal embedding such that $P=\Phi (x_q)$ and let $\mathfrak a\in
\I (R)$. There exists $\beta \in \cO_{D, N}$, $\n(\beta )>0$, such
that $q(\mathfrak a)\cO_{D, N} = \beta \cO_{D, N}$ and for any such
$\beta$ we have
$$
P^{\sigma _{\mathfrak a}} = \Phi (\beta^{-1} x_q).
$$

\end{enumerate}

\end{teor}

$\\ $ By part (1) of the above Theorem, we know that $H_R$ is an
extension of $Q(P)$ of degree at most $2$. The determination of this
subfield of $H_R$ is the main result of this section, which is
contained in Theorem \ref{mainCM}. In order to obtain this, we need
some previous lemmas.

\begin{lema}\label{calculd}
Let $m\mid D\cdot N$ and $P\in \mathrm{CM}(R)$. Then
$\om_m(P)=P^{\sigma}$ for some $\sigma \in \Gal(H_R/K)$ if and only
if $m\mid \frac{D N}{D(R) N(R)}$. If this is the case,
$\sigma={\sigma_{\mathfrak b}}$ for the ideal $\mathfrak b\in \I (R)$
such that $\Norm_{K/\Q} (\mathfrak b) = m$.
\end{lema}

\noindent {\bf Proof. } If $p\mid D(R) N(R)$ then $\om_p$ switches the local orientation of $P$ at $p$, whereas
any $\sigma \in \Gal(H_R/K)$ preserves it and hence $\om_p(P)\not \in \Gal(H_R/K)\cdot P$.

Let $p\mid \frac{D N}{D(R) N(R)}$. By Proposition \ref{og}, $p$ is a
ramified prime of $R$ which does not divide the conductor $f$ and
thus there exists  $\mathfrak p\in \I (R)$ such that $\Norm_{K/\Q
}(\mathfrak p) = p$. Let $\beta \in \cO_{D, N}$, $\n(\beta )>0$, be
such that $q(\mathfrak p)\cO_{D, N} = \beta \cO_{D, N}$. Since
$\Norm_{K/\Q }(\mathfrak p) = p$ we have $\n (\beta )=p$.

Let $P=\Phi(x_q)$ for some optimal embedding $q:R\hookrightarrow \cO_{D, N}$. By Theorem \ref{Shi},
we have $P^{\sigma _{\mathfrak p}}=\Phi (\beta^{-1}x_q)$. Let us show that $\om_p(P) = \Phi (\beta x_q) = \Phi (\beta^{-1}x_q)$ as well.
In order to prove this, it follows from the definition and properties of the Atkin-Lehner group that it suffices to show that
$\cO_{D, N}\cdot \beta $ is a two-sided ideal.

For $p\mid \frac{D}{D(R)}$ this is immediate, as there exists a
unique ideal of norm $p$ in $\cO _{D, N}$, which is two-sided. Assume now
that $p\mid \frac{N}{N(R)}$. The question is local and for all
primes $\ell\ne p$, $\cO_{D, N}\otimes \Z_{\ell }\cdot \beta =
\cO_{D, N}\otimes \Z_{\ell }$. Locally at $p$, the Eichler order
$\cO_{D, N}\otimes \Z_p$ can be described as $\cO_{D, N}\otimes \Z_p
= \{ \begin{pmatrix} a&b\\c\, p&d
\end{pmatrix}: a, b, c, d\in \Z_p\}$. Since $\mathfrak p^2 =(p)$ it follows that $\cO_{D, N}\otimes \Z_p\cdot \beta^2 = p\cdot \cO_{D, N}\otimes \Z_p$ and it is elementary to check that $\beta \in (\cO_{D, N}\otimes \Z_p)^*\cdot \begin{pmatrix}
0&1\\p&0
\end{pmatrix}$ and hence normalizes $\cO_{D, N}$. $\Box $

We next describe the action of complex conjugation on points of
complex multiplication on $X_0(D, N)$. This shall serve us in the
next section to exhibit rational points on the quotient curves
$X_0(D, N)/\langle \om _m\rangle $.

\begin{lema}\label{CC}
Let $P\mapsto \bar {P}$ denote the complex conjugation on $X_0(D,
N)(\qbar )$. Then

\begin{enumerate}

\item $\overline{\mathrm{CM}(R)} = \mathrm{CM}(R)$.

\item For every $P\in \mathrm{CM}(R)$ there exist  $m\mid D(R)\, N(R)$ and
$ \mathfrak a \in \I (R)$ such that $$\bar P= \om_{m} (P^{\sigma
_{\mathfrak a}})\,. $$
Moreover, the integer   $m $ and the class $[\ga]\in
\Pic(R)/\Pic(R)^2$ do not depend on $P$. More precisely, $m=D(R)\,
N^*(R)$,   $\mathfrak a$ satisfies
$$B_D\simeq (\frac{-s, D(R) N^*(R) \mbox{N}_{K/\Q }(\mathfrak a)}{\Q })\,,$$
and for all $\mathfrak b\in[\mathfrak a]$ there exists $Q\in  \Gal(H_R/K)\cdot P$ such that $\bar Q=\om_m(Q^{\sigma_{\mathfrak b}})$.

\end{enumerate}

\end{lema}

\noindent {\bf Proof.} Let $P=\Phi(x_q)$ for some normalized optimal
embedding $q:R\hookrightarrow \cO_{D, N}$. Let $\epsilon \in \cO_{D,
N}$ be of reduced norm $\n(\epsilon )=-1$. By \cite{Sh2},
\cite{Ogg1}, $\bar P = \Phi (\epsilon \bar x_q)$.

$(1)$ We have $\bar P = \Phi (\epsilon \bar x_q) = \Phi (x_{q\cdot
\epsilon^{-1}})$. The embedding $a\mapsto \epsilon q(\bar a)
\epsilon^{-1}$ is a normalized optimal embedding of $R$ into
$\cO_{D, N}$ and hence $\bar P$ is also a point of complex
multiplication by $R$.

(2) Let $P\in \mathrm{CM}(R)$. By Proposition \ref{og}, $\bar P = \om_{m_P} (P^{\sigma _{\mathfrak a_P}})$ for some $m_P\mid D(R)\, N(R)$ and $\mathfrak a_P\in \I (R)$. For a fixed $P$,
the integer $m_P$ and the class of the ideal $\mathfrak a_P$ in $\Pic (R)$ are unique, but may
depend on $P$.

Let $Q\in \mathrm{CM}(R)$ be any other point with complex
multiplication by $R$. By Proposition \ref{og}, $Q=\om_{d}
(P^{\sigma _{\mathfrak b}})$ for some $d\mid D(R)\,N(R)$ and $\mathfrak b\in \I (R)$. Thus $\bar Q = \om_{d}\overline{(P^{\sigma _{\mathfrak b}})} = \om_{d}(\bar{P})^{\sigma^{-1}_{\mathfrak b}} =
\om_{d}(\om_{m_P} (P^{\sigma _{\mathfrak
a_P}}))^{\sigma^{-1}_{\mathfrak b}} = \om_{m_P}(Q^{\sigma_{\mathfrak a \mathfrak b^{-2}}})$. Hence, $m_P=m_Q$ and  $[\mathfrak a_P ]= [\mathfrak a_Q]$ in $\Pic(R)/\Pic(R)^2$ for
any two points $P, Q\in \mathrm{CM}(R)$.

Let $m=m_P$ and $\mathfrak a = \mathfrak a_P$.
Changing, if  necessary, $\mathfrak a$ by an ideal in the same class in $\Pic (R)$, we can assume  that $q(\mathfrak a^{-1})\subset \cO_{D, N}$ and $(\Norm_{K/\Q }(\mathfrak
a^{-1}),D N)=1$. For this choice,  let $\beta \in \cO_{D, N}$ be such that
$q(\mathfrak a^{-1})\cO_{D, N} = \beta \cO_{D, N}$.
Let $i=q(\sqrt{-s})\in \cO_{D, N}$.  We first show the following

$\\ $ {\bf Claim.} {\em There exists $j\in \cO_{D, N}$ such that $j^2=m
\Norm_{K/\Q }({\mathfrak a^{-1}})$, $i j=-j i$.}

We have  that $\bar P = \om_{m} (P^{\sigma _{\mathfrak a}})$. Since
$\bar P = \Phi (\epsilon \bar x_q)$, by Theorem \ref{Shi}
$$\Phi (\epsilon \bar x_q) = \Phi (\gamma \beta x_q)$$
for some $\gamma \in \cO_{D, N}$, $\n (\gamma )=m$. Thus $\alpha
\epsilon \bar x_q = \gamma \beta x_q$ for some $\alpha \in \cO^1_{D,
N}$. If we write $\eta = \alpha \epsilon $, this reads
$$\eta \bar x_q = \gamma \beta x_q, \mbox{ with }\n(\eta
)=-1.$$

Let $j =\eta ^{-1}\gamma \beta \in \cO_{D, N}$. As it is checked,
$j(\cH )=\bar \cH $ and $j(\bar \cH )=\cH$; with $j(\bar x_q)=x_q$
and $j(\bar x_q)=x_q$. Hence $j$ has no fixed points on $\cH \cup
\bar \cH$ whereas $j^2$ fixes $x_q$ and $\bar x_q$. By
\cite[Proposition 1.2]{Sh2} we obtain that $j^2\in \Q ^*$. Since
$\n(j)=-m \mbox{N}_{K/\Q }({\mathfrak a^{-1}})\in \Z $, $j^2 = m
\mbox{N}_{K/\Q }({\mathfrak a^{-1}})$.

Moreover, the single fixed point of $j q j^{-1}:R\hookrightarrow
\cO_{D, N}$ on $\cH $ is $j(\bar x_q) = x_q$. Hence, either $j
q(\sqrt{-s})j^{-1} = q(\sqrt{-s})$ or $-q(\sqrt{-s})$. In the first
case, it implies that $j\in q(K)$ and this can not be possible,
since $\n(j)<0$. Hence, $i j = -j i $. This proves our claim.

Let us now show that $m=D(R)\,N^*(R)$. By the discussion
following Proposition \ref{og}, it suffices to show that for $p\mid
D(R) N(R)$, we have $o_p(P)=-o_p(\bar P)$ if and only if $p\mid D(R)
N^*(R)$. In other words, for $p\mid D(R) N(R)$, there exists $\alpha
\in (\cO_{D, N}\otimes \Z_p)^*$ such that $\bar q = \alpha q
\alpha^{-1}$ if and only if $p\mid (N, f)$.

Let $p\mid D(R)$. Then $B\otimes \Q_p = q(K\otimes \Q_p)\oplus
q(K\otimes \Q_p)\cdot \pi $, where $\pi^2=p$ and $a\cdot \pi = \pi
\cdot \bar a$ for any $a\in K$. For any $\alpha \in B_p$, we have
$\alpha q(\sqrt{-s})\alpha^{-1} = -q(\sqrt{-s})$ if and only if
$\alpha \in q(K\otimes \Q_p)\cdot \pi $, but none of these elements
is an integral unit. Thus $o_p(P)=-o_p(\bar P)$.

Let $p\mid N(R)=N^*(R)\cdot (N, f)$. Then $B_p\simeq \M_2(\Q_p)$ and
we may assume that $\cO_{D, N}\otimes \Z_p = \{
\smaller{\begin{pmatrix} a & b\\ p \,c & d\end{pmatrix}}: a, b, c,
d\in \Z_p\} \subset \M_2(\Q_p)$. Assume now that $p\mid N^*(R)$.
Then $q: K\otimes \Q_p \simeq \Q_p\times \Q_p \ra \M_2(\Q_p)$, $(a,
d)\mapsto \smaller{\begin{pmatrix} a&0\\0&d\end{pmatrix}}$. Since
$\bar q (a, d) = q(d, a)$, we have $\alpha q\alpha^{-1} = \bar q$ if
and only if $\alpha =\smaller{\begin{pmatrix} 0 & b \\ p\, c & 0
\end{pmatrix}}$. None of these elements are units of $\cO_{D,
N}\otimes \Z_p$, hence $o_p(P)=-o_p(\bar P)$.

On the other hand, if $p\mid (N, f)$ we have (unless $p\ne 2$ and $s
\equiv 3$ mod $4$) $R\otimes \Z_p \simeq \Z_p + p^k\Z_p \sqrt{-s}$,
where $p^k\parallel f$, $k\ge 1$. We can write $q:R\otimes \Z_p\ra
\cO_{D, N}\otimes \Z_p$,
$q(p^k\sqrt{-s})\mapsto \smaller{\begin{pmatrix} 0 & 1 \\
-p^{2k} s & 0
\end{pmatrix}}$. Since $\smaller{\begin{pmatrix} 0 & -1 \\ p^{2k} s & 0
\end{pmatrix} = \begin{pmatrix} 1 & 0 \\ 0 & -1
\end{pmatrix}\begin{pmatrix} 0 & 1 \\ -p^{2k} s & 0
\end{pmatrix}\begin{pmatrix} 1 & 0 \\ 0 & -1
\end{pmatrix}^{-1}}$, we deduce that $o_p(P)=o_p(\bar P)$.

Similarly, if $p=2\mid(N, f)$ and $s\equiv 3$ mod $4$, we have
$R\otimes \Z_2 \simeq \Z_2 + 2^{k-1}\Z_2 \sqrt{-s}$, where
$2^k\parallel f$, $k\ge 1$, and $q(2^{k-1}\sqrt{-s})\mapsto \smaller{\begin{pmatrix} 2^{k-1} & -2^{k-2}(s+1) \\
2^k & -2^{k-1}
\end{pmatrix}}$. As before, conjugating by $\smaller{\begin{pmatrix} 1 & 0 \\ 0 & -1
\end{pmatrix}}$ we obtain that $o_p(P)=o_p(\bar P)$. Summing up, we
have shown that $m=D(R)\,N^*(R)$.

For the ideal $\mathfrak a$, we have
$B_D=\Q (\,i, j\,) = (\frac{-s, m \mbox{N}_{K/\Q }(\mathfrak a^{-1})}{\Q
})\simeq(\frac{-s, m \mbox{N}_{K/\Q }(\mathfrak a)}{\Q })$.

Finally, let $\tau = \sigma_ {\mathfrak a}\cdot \sigma_0^{-2}$ with $\sigma_0\in\Gal(H_R/K)$. Then, $\bar Q =\om_{m} (Q^{\tau })$,
where $Q=P^{\sigma_0}$.
\hfill $\Box $

\begin{remark}\label{uniqueness}
If  $\mathfrak a\in\I (R)$ satisfies $B_D\simeq(\frac{-s, m
\operatorname{N}_{K/\Q }(\mathfrak a)}{\Q })$, then any other ideal
in the class of $[\mathfrak a]\in\Pic(R)/\Pic (R)^2$ also satisfies
this isomorphism. In general, the converse is not true, but if $H_R$
is the Hilbert class field of $K$, then $[\mathfrak
a]\in\Pic(R)/\Pic (R)^2$ is uniquely determined. Indeed, the
isomorphism $(\frac{-s, m \operatorname{N}_{K/\Q }(\mathfrak a)}{\Q
})\simeq (\frac{-s, m \operatorname{N}_{K/\Q }(\mathfrak b)}{\Q })$
implies that $\operatorname{N}_{K/\Q}( \mathfrak a\cdot \mathfrak b
^{-1})=\operatorname{Norm}_{K/\Q}(\alpha)$ for some $\alpha \in
K^*$. Then the value $ \operatorname{Norm}_{K/\Q}(\alpha) \pmod
{\operatorname{disc} (R)}$ is represented by a quadratic form of the
principal genus of discriminant $\operatorname{disc}(R)$ and by
genus theory $[ \mathfrak a]=[ \mathfrak b]$.
\end{remark}

As a consequence of Theorem \ref{Shi} and Lemma \ref{CC}, we obtain
the following result.

\begin{teor}\label{mainCM}

Let $P\in \CM(R)$. We have that

\begin{enumerate}

\item If $D(R) N^*(R) \ne 1$ then $\Q(P) = H_R$.

\item If $D(R) N^*(R) = 1$ then $[H_R:\Q(P) ]=2$ and $\Q(P) \subset H_R$ is the
subfield fixed by $\sigma = c\cdot \sigma_{\mathfrak a}\in \Gal
(H_R/\Q )$ for some  $\mathfrak
a\in\I(R)$ such that $B_D\simeq (\frac{-s,
\operatorname{N}_{K/\Q }(\mathfrak a)}{\Q })$, where $c$ denotes the complex conjugation.

\end{enumerate}

\end{teor}

\noindent {\bf Proof.}  Let $P\in \mathrm{CM}(R)$. By Lemma \ref{CC}
(2),
$$
\Gal (H_R/\Q )\cdot P = \left(\Gal (H_R/K)\cdot P\right) \cup
\left(\Gal (H_R/K)\cdot \om_{D(R)N^*(R)}(P)\right).
$$
Assume first that $D(R)N^*(R)\ne 1$. Then $$\left(\Gal (H_R/K)\cdot
P\right) \cap \left(\Gal (H_R/K)\cdot \om_{D(R)N^*(R)}(P)\right) =
\emptyset $$ because the action of $W(R)\times \Pic(R)$ on
$\mathrm{CM}(R)$ is free by Proposition \ref{og}. Hence $\sharp \Gal
(H_R/\Q )\cdot P = 2 h(R)=[H_R:\Q ]$. Since $K(P)\subseteq H_R$ by
Theorem \ref{Shi}, it follows that $\Q (P)=H_R$.

Assume now that $D(R)N^*(R)=1$. By Lemma \ref{CC} (2),  $\sharp \Gal (H_R/\Q )\cdot P =
 h(R)$ and $\Q (P)\subset H_R$ must be a subfield of
index $2$ of $H_R$. Again, by Lemma \ref{CC} (2), $\Q (P)$ is the subfield
of $H_R$ fixed by $c\cdot \sigma_{\mathfrak a}\in \Gal (H_R/\Q )$ for some  $\mathfrak a\in\I (R)$ such that $B_D\simeq (\frac{-s, \mbox{N}_{K/\Q
}(\mathfrak a)}{\Q })$. \hfill $\Box$

\begin{cor}\label{d}
Assume that either $h(R)=1$ or $h(R)=2$ and $D(R) N^*(R) = 1$. Then
$X_0(D,N)$ admits rational points over some imaginary quadratic field.
\end{cor}

Let $m\vert D\cdot N$ and let $\pi_m:X_0(D,N)\ra X_0(D,N)/\langle
\om_m\rangle $ denote the natural projection map. We say that a
point $Q\in X_0(D,N)/\langle \om_m\rangle (\qbar )$ is a CM point if
$\pi_m^{-1}(Q)$ is a pair of CM points on $X_0(D,N)$.

\begin{cor}\label{AL}
Let $P\in \CM(R)\subset X_0(D,N)(\qbar )$ and $Q=\pi_m(P)$ for some $m\mid D N$. Set $m_r= \gcd (m, \frac{D\,N}{D(R)\,N(R)})=\gcd(m,\operatorname{disc}(R)/\gcd (N,f))$ and
 let $\mathfrak b$ be the invertible ideal of $R$
such that $\Norm_{K/\Q }(\mathfrak b)=m_r$.

\begin{enumerate}

\item Assume $D(R)N^*(R)\ne 1$. Then $\Q(Q)$ is
$$\left\{\begin{array}{ll} H_R^{\sigma_{\mathfrak b}} & \text{ if } m/m_r=1, \\
H_R^{\sigma_{\mathfrak b \mathfrak a}\cdot c} \text{ for some
}\mathfrak a \text{ such that } B_D\simeq (\frac{-s, \mbox{N}_{K/\Q
}(\mathfrak a)}{\Q })& \text{ if } m/m_r=D(R) N^*(R),
\\
H_R& \text{ otherwise.}  \end{array}\right.$$

\item Assume $D(R)N^*(R)=1$. Then

$$
\Q(Q)=\left\{\begin{array}{ll} H_R^{\langle c\cdot\sigma_{\mathfrak
a},
\sigma_{\mathfrak b}\rangle }& \text{ if } m/m_r=1,  \\
H_R^{c\cdot\sigma_{\mathfrak a}} & \text{ otherwise,}
\end{array}\right.$$
for some $\mathfrak a\in\Pic (R)$ such that $B_D\simeq (\frac{-s,
\mbox{N}_{K/\Q }(\mathfrak a)}{\Q })$.
\end{enumerate}

\end{cor}

\end{document}